\documentclass[12pt]{article}

\usepackage{amssymb}
\usepackage{amsmath}
\usepackage{amscd}
\usepackage{latexsym}
\usepackage{graphicx}

\usepackage{cite}

\newcommand{\be}{\begin{equation}}
\newcommand{\ee}{\end{equation}}
\newcommand{\Dlt}{\Delta}

\newcommand{\al}{\alpha}
\newcommand{\ra}{\rightarrow}

\newcommand{\gm}{\gamma}

\begin{document}

\begin{center}

{\Large{\bf Self-similarly corrected Pad\'{e} approximants for the indeterminate problem} \\ [5mm]

S. Gluzman and V.I. Yukalov }  \\ [3mm]

{\it Bogolubov Laboratory of Theoretical Physics, \\
Joint Institute for Nuclear Research, Dubna 141980, Russia} 

\end{center}

\vskip 5cm

\begin{abstract}

A method is suggested for treating the well-known deficiency in the use of
Pad\'e approximants that are well suited for approximating rational functions,
but confront problems in approximating irrational functions. We develop the 
approach of self-similarly corrected Pad\'e approximants, making it possible 
to essentially increase the class of functions treated by these approximants. 
The method works well even in those cases, where the standard Pad\'e approximants 
are not applicable, resulting in divergent sequences. Numerical convergence of
our method is demonstrated by several physical examples. 

\end{abstract}

\vskip 2cm
{\parindent=0pt
{\bf PACS numbers}: 02.30.Lt (Sequences, series, and summability); 
02.30.Mv (Approximations and expansions); 02.90.+p (Other topics in mathematical 
methods in physics) }

\newpage

\section{Introduction}

In many branches of physics and applied mathematics, there exists the standard 
problem, when a sought function is given as an expansion in powers of an 
asymptotically small variable, but it is necessary to define this function 
in the domain of finite, or large, or even infinite values of the variable. 
Then some extrapolation techniques are required.

Probably, the simplest and the most often used method of extrapolation is
by means of Pad\'{e} approximants. An approximant $P_{n/m}(x)$ is the ratio 
of two polynomials $P_{n}(x)$ and $P_{m}(x)$ of the order $n$ and $m$, 
respectively. The coefficients of the Pad\'{e} approximant are defined
trough the coefficients of the given power series \cite{b1,Suetin}, from the
requirement of the asymptotic equivalence to the given series of the sought 
function $f(x)$.

Pad\'{e} approximants locally are the best {\it rational} approximations of
power series. Their poles determine singular points of the approximated
functions \cite{b1,Suetin}. Calculations with Pad\'{e} approximants are 
straightforward and can be performed with 
$Mathematica^{\textsuperscript{\textregistered}}$.

From the emerging sequences of approximants, one has to select those that are
holomorphic functions. Not all approximants generated by the procedure are
holomorphic. The holomorphy of diagonal Pad\'{e} approximants in a given domain
implies their uniform convergence inside this domain \cite{gonch}.

However, the use of the Pad\'{e} approximants confronts the well known 
{\it indeterminate problem}, when these approximants are not applicable. The 
meaning of this problem is as follows.

Suppose we are considering a non-negative function with the asymptotic behavior 
at infinity $f(x)\simeq A x^{s}$, with a known {\it critical index} $s$ and 
a given expansion at small $x$,
\begin{equation}
\label{small}
f(x)=a_{0}+\sum _{n=1}^N a_n x^n+O(x^{N+1}) \; .
\end{equation}
Here $N$ is an integer and $N \geqslant 1$. The standard suggestion for the
solution of the problem of reconstructing the amplitude $A$, is based on the
following \cite{Baker_3,bender}.

Suppose we need to calculate the critical amplitude $A$. To this end, let us
apply a transformation of the original series $f(x)$, so that to obtain the
transformed series $T(x) = f(x)^{-1/s}$, in order to get rid of the $x^{s}$
behavior at infinity. Applying the technique of diagonal Pad\'{e} approximants
in terms of $x T_{n/(n+1)}(x)$, one can readily obtain the sequence of approximations
${A_{n}}$ for the critical amplitude $A$ as
\begin{equation}
\label{sug1}
A_{n}=\lim_{x\to \infty } (x T_{n/(n + 1)}(x))^{-s} \; .
\end{equation}

Of course, not always Pad\'{e} approximants are convergent. The uniform 
convergence of diagonal sequences of Pad\'{e} approximants has been
established, e.g., for the Stieltjes-Markov-Hamburger series \cite{b1,lub}.
The Stieltjes moment problem can possess a unique solution or multiple solutions.
The latter is called the {\it indeterminate problem}. The existence of a unique
solution or multiple solutions depends on the behavior of the moments. This is in 
contrast to the problem of moments for a finite interval \cite{kr}, which is 
uniquely defined, if the solution exists \cite{Wall}. The original work of 
Stieltjes has been clearly explained in \cite{st,st1}.  

From the theory of $S$-fractions one deduces a {\it divergence by oscillations} 
of Pad\'{e} approximants  in the indeterminate case \cite{Wall}. The phenomenon
of spurious poles, when approximants can have poles which are not related to the
underlying function, can also lead to divergence of the Pad\'{e} sequence
\cite{b1,lub,Wall}.

Generally, the form of the solution in the indeterminate case is covered by the
Nevanlinna theorem, expressing non-negative functions through a mixture of
rational and non-rational contributions\cite{kr,ad,adam,mal}.
The former can be interpreted as the part covered by four entire functions related
to the construction of Pad\'{e} approximants, while the latter is only a
non-rational function from the Nevanlinna-class \cite{berg}. Thus, in the
indeterminate case, Pad\'{e}  approximants are not sufficient. Below, we suggest 
an approximation scheme which can express the irrational part explicitly and 
describes the rational part with Pad\'{e} approximants.

In the present paper, we consider approximations for single-valued real functions, 
because of which we keep in mind the standard Pad\'{e} approximants. We will 
not consider complex multivalued functions, or sets of functions, which are all 
given as Taylor series around the same point, when we would need to find a 
simultaneous rational approximation to such a vector of several functions, when 
it would be possible to employ the Pad\'{e} - Hermite approximation 
\cite{b1,Assche_36} yielding the so-called algebraic approximants 
\cite{b1,Assche_36,Sergeev_37}. Note that, instead of considering simultaneously 
a set of functions, it is admissible to treat each of them separately.

\section{Self-similarly corrected Pad\'{e} approximants}

The main idea of the method we suggest is to separate out an initial approximation
that could treat the irrational part of the sought function, so that the higher 
approximations, constructed above the initial approximation, could be represented 
by Pad\'{e} approximants. The form of the initial approximation can be chosen by
resorting to self-similar approximation theory \cite{Y_1,Y_2} resulting in root
approximants \cite{g1,g2} or factor approximants \cite{f1,f2}.

Assume that we start with a self-similar approximation $K(x)$, ensuring the 
correct critical index $s$, while all other additional parameters are to be 
obtained by asymptotically matching the approximant with the truncated series 
for $f(x)$. As a starting approximation $K(x)$, we can employ iterated root 
approximants \cite{g1,g2} or factor approximants \cite{f1,f2}. Then the 
initial approximation for the amplitude is 
$$
A_{0} = \lim_{x \to \infty} (K(x) x^{-s}) \; .
$$
To define the higher-order approximations, we can resort to the method of  
{\it corrected approximants} \cite{g1, g2}. For this purpose, to find out 
corrections to the critical amplitude, we divide the original series
(\ref{small}) for $f(x)$ by $K(x)$, and denote the new series as
$$
G(x) = \frac{f(x)}{K(x)} \; .
$$
Assuming that the irrational part has already been included in $K(x)$, we 
now can invoke rational approximants. Thus, we finally construct a sequence 
of diagonal Pad\'{e} approximants asymptotically equivalent to $G(x)$, so that 
the approximate amplitudes are expressed by the formula
\begin{equation}
\label{main}
A_{n}=A_{0} \lim_{x\to \infty} [G_{n/n}(x)] \; .
\end{equation}
The role of the starting approximation $K(x)$ is crucial. It is supposed not
only to approximate the irrational part of the sought function, but also to
ensure the convergence of Pad\'{e} approximants for the rational part.

The starting approximation can be viewed as a {\it control function}, constructed
so that to ensure the convergence of the sequence of corrected Pad\'{e}
approximants, even if the standard Pad\'{e} scheme diverges. Our the most important
suggestion is not only in promoting the idea of controlled Pad\'{e} sequences, but
also in conjecturing that the control function should be chosen among the
low-parametric subset of irrational self-similar approximants.

In the examples to be treated below, we assume that the critical index at infinity 
is available. This information, as is known, makes it possible to improve the 
accuracy of approximants \cite{Baker_3,g2,Y_39,Fernandez_38}.

When only a single "critical index" $s$ is known in advance, we are limited
to the iterated roots and factor approximants, noting that self-similar
approximants can be uniquely defined. In all examples studied below, the
control function is indeed found easily. More general self-similar roots and
continued roots \cite{cross,cont,Y_40}, can be applied to more general situations, 
when, in addition to $s$, more detailed corrections to the leading
scaling behavior are available. 

Note that Pad\'{e} approximants by themselves cannot serve as control functions,
since the standard Pad\'{e} scheme cannot be controlled \cite{rep,rep1}.
It is also useful to mention that the suggested approach can be generalized 
to the case of complex functions, for which self-similar approximation theory
is applicable, with the use of complex control functions \cite{Y_39}. However,
here we limit ourselves by single-valued real functions.   

In this way, the corrected approximants are applicable in all cases, including 
those where the standard Pad\'{e} approximants diverge. On the other hand, when 
the standard scheme is convergent, the corrected scheme also converges, exhibiting 
either the similar behavior of approximants, or often accelerating convergence. 
These cases are illustrated in the following sections.

Let us emphasize that our main aim in the present paper is to extend the 
applicability of Pad\'{e} approximants. Of course, we could employ one of the 
variants of self-similar approximation theory \cite{Y_1,Y_2,g1,g2,f1,f2,cross,cont}
for the whole problem, not invoking Pad\'{e} approximants at all. However, the 
techniques of Pad\'{e} approximants are well developed and there are standard 
programs for defining their coefficients. It would therefore be tempting to 
slightly modify the method of Pad\'{e} approximants in order to be able to resort
to the standard techniques of calculating their coefficients, at the same time 
extending their applicability to the cases where the standard Pad\'{e} approximants 
are not applicable. Such an extension is a principal novelty of the present paper,
not considered in our previous publications.

\section{Self-similar roots as control functions}

Choosing an irrational part $K(x)$, playing the role of a control function, it
is convenient to invoke the self-similar root approximants \cite{g1,g2,Y_40},
which we shall employ below. To make this paper self-consistent, we briefly  
delineate the way of defining the self-similar root approximants. A more detailed
description can be found in the earlier publications \cite{g1,g2,Y_40}. 

Suppose, the sought function $f(x)$ is given as an expansion in powers of a small 
variable $x \ra 0$, so that the $k$-th order series reads as
\be
\label{S1}
 f_k(x) = f_0(x) \left ( 1 + \sum_{n=1}^k a_n x^n \right ) \;  ,
\ee
where 
\be
\label{S2}
f_0(x) = A x^\al \qquad ( A \neq 0 ) 
\ee
is a known zero-order term. And assume that we need to extrapolate this series
to the region of finite variables, including the limit of the large variable, 
where the sought function behaves as
\be
\label{S3}  
 f(x) \simeq B x^s \qquad ( x \ra \infty )  ,
\ee
with an unknown amplitude $B$. Then the $k$-th order self-similar root approximate 
$R_k^*(x)$ is given \cite{g1,g2,Y_40} by the formula
\be
\label{S4}
 \frac{R_k^*(x)}{f_0(x)} = \left ( \left ( \ldots ( 1 + A_1 x)^{n_1} +
A_2 x^2 \right )^{n_2} + \ldots + A_k x^k \right )^{n_k} \; ,
\ee
in which the first $k-1$ powers are
\be
\label{S5}
n_j = \frac{j+1}{j} \qquad ( j = 1,2,\ldots , k-1) \;   ,
\ee
while the power $n_k$ is
\be
\label{S6}
 n_k = \frac{s-\al}{k} \qquad ( k = 1,2,\ldots ) \;  .
\ee
Here all parameters $A_i$ are uniquely defined by the accuracy-through-order 
procedure, that is, by expanding this root approximant in powers of $x$ and 
equating the similar terms of this expansion and of series (\ref{S1}).  

The root approximant (\ref{S4}) at large variables behaves as
\be
\label{S7}
 R_k^*(x) \simeq B_k x^s \qquad ( x \ra \infty ) \;  ,
\ee
with the amplitude
\be
\label{S8}
 B_k = A \left ( \left ( \ldots ( A_1^2 + A_2)^{3/2} + A_3 \right )^{4/3} +
\ldots + A_k \right )^{(s-\al)/k} \;  .
\ee

In order to separate the irrational part of the sought function, it is sufficient
to employ a low-order root approximant, usually, accepting that of second order.
Below, we illustrate the method by several examples. First, we consider rather
simple functions, for which exact forms are known, and respectively, for which the
accuracy of the method can be explicitly estimated. Then we treat more complicated
cases related to interesting physical problems.

\section{Convergence and accelerated convergence}

\subsection{Convergence: Mittag-Leffler function}

Consider the Mittag-Leffler function appearing in fractional evolution 
processes \cite{mittag},
\begin{equation}
\label{mittag}
F(x)= \rm{erfc}(x) \exp \left(x^2\right) \; .
\end{equation}

For small $x$ in the lowest orders, one has
\begin{equation}
F(x)\simeq 1-\frac{2 x}{\sqrt{\pi }}+x^2-\frac{4 x^3}{3 \sqrt{\pi }}+\frac{x^4}{2}+O(x^5) \; ,
\end{equation}
and for large $x$, 
$$
F(x) \simeq A x^{-1} \qquad (x\rightarrow \infty) \; ,
$$
where $A = 1/\sqrt{\pi }$

Let us take as a control function the second-order self-similar root approximant
\begin{equation}
K(x)=R_{2}^{*}(x) = 
\left[ \left( \frac{2 x}{\sqrt{\pi }}+1\right)^2-\frac{2 (\pi -4) x^2}{\pi} \right]^{-1/2} \; .
\end{equation}

\begin{figure}
\begin{center}
\includegraphics[width=10cm]{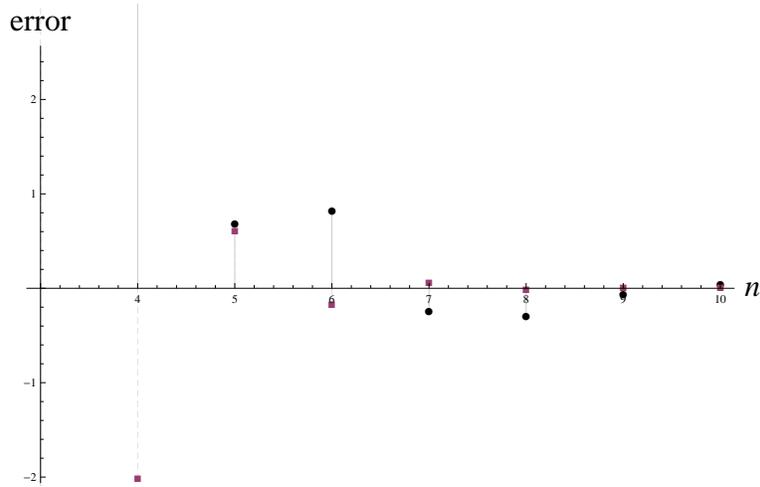}
\end{center}
\caption{
Percentage errors for the large-variable amplitude of the Mittag-Leffler function, 
with the increasing approximation order $n$: the errors for the corrected Pad\'{e} 
approximants are shown with circles, and for the standard Pad\'{e} scheme, with 
squares.
\label{fig:mitt}}
\end{figure}

There is no significant difference in the behavior of approximations for the 
amplitude $A$ for both schemes of applying Pad\'{e} approximants. Good convergence 
is achieved as is shown in Fig. \ref{fig:mitt}.

\subsection{Accelerated convergence: quartic oscillator}

Let us consider the dimensionless ground state energy $e(g)$ of the quantum
one-dimensional quartic anharmonic oscillator \cite{osc1}. Here $g$ stands for
the dimensionless coupling constant. The asymptotic expansion of $e(g)$ in the
weak-coupling limit, when $g\rightarrow 0$, is
\begin{equation}
\label{30}
e(g)\simeq a_0+a_1g+a_2g^2+...+a_{18}g^{18}+... \; ,
\end{equation}
where a few starting coefficients are
$$
a_0 = \frac{1}{2}, \quad a_1 = \frac{3}{4}, \quad a_2 = - \frac{21}{8}, \quad 
a_3 = \frac{333}{16},   \quad a_4 = - \frac{30885}{128} \; .
$$
At a large coupling parameter $g \rightarrow \infty$, the ground-state energy $e(g)$
diverges as $A g^{1/3}$, with $A = 0.667986$.

Applying the standard Pad\'{e} approximants, one obtains the amplitudes
$$
A_{2}=0.759147, \quad A_{3}=0.734081, \quad A_{4}=0.720699, \quad A_{5}=0.712286, 
$$
$$
A_{6}=0.706466, \quad A_{7}=0.702176, \quad A_{8}=0.698869, \quad A_{9}=0.696173 \; . 
$$
The last approximant of the sequence has the error of $4.21967\%$.

The control function for the method of corrected Pad\'{e} approximants is the 
second-order iterated root,
\begin{equation}
K(g)=R_{2}^{*}(g)=\frac{1}{2} \left( \left(\frac{9 g}{2}+1\right)^2-18 g^2 \right)^{1/6} \; ,
\end{equation} 
yielding $A_{0}=0.572357$. 

The results of calculations, according to the corrected Pad\'{e} approximants
methodology, are as follows: 
$$
A_{1}=A_{2}=A_{0}, \quad A_{3}=0.587104, \quad A_{4}=0.63279, \quad A_{5}=0.655086, 
$$
$$
A_{6}=0.660334, \quad A_{7}=0.661945, \quad A_{8}=0.663225, \quad A_{9}=0.665346 \; .
$$
The last approximant of the sequence gives an error of only $-0.3952\%$, by the order 
of magnitude better than the standard approach.

\section{Treating divergence by oscillations}

\subsection{Correlation function}

In the course of calculating temporal correlation functions \cite{mori}, one often 
meets the functions of the following structure 
\begin{equation}
\label{mor}
f(x)=\sqrt{x^2+4}-x \; .
\end{equation}

At small $x$, the expansion is 
\begin{equation}
f(x)\simeq 2-x+\frac{x^2}{4}-\frac{x^4}{64}+O(x^6) \; .
\end{equation}

For large $x$, the following asymptotic behavior holds:
\begin{equation}
f(x)\simeq 2 x^{-1}  + O(x^{-3})\; .
\end{equation}

\begin{figure}
\begin{center}
\includegraphics[width=10cm]{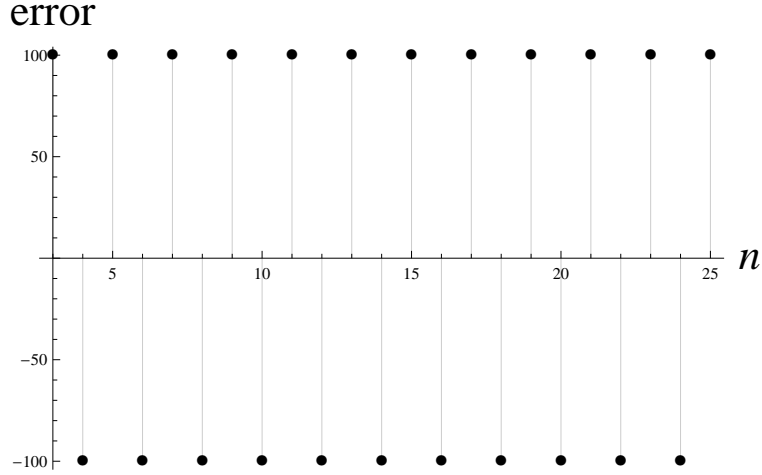}
\end{center}
\caption{
Percentage errors for the amplitude, with increasing approximation order, 
obtained with Pad\'{e} approximants for function (\ref{mor}).
\label{fig:mori1}}
\end{figure}

Applying the standard Pad\'{e} scheme brings a divergent by oscillation
solution for the amplitude, bounded by the values $0$ and $4$, as is shown
in Fig. \ref{fig:mori1}.

\begin{figure}
\begin{center}
\includegraphics[width=10cm]{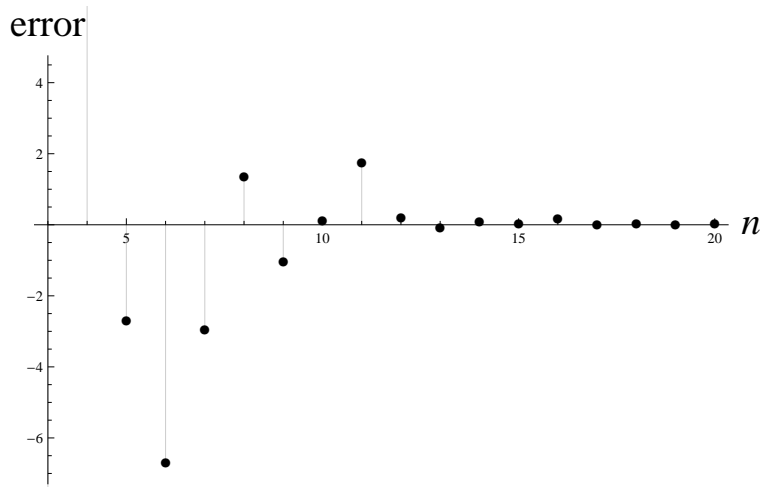}
\end{center}
\caption{
Percentage errors for the amplitude, with increasing approximation order, 
obtained with the corrected Pad\'{e} approximants for function (\ref{mor}).
\label{fig:mori2}}
\end{figure}

But employing the corrected scheme with the control function
\begin{equation}
K(x)=R_{2}^{*}(x)=\frac{2}{\sqrt{x^2/2+x+1}} 
\end{equation}
leads to a convergent result for the amplitude, as is shown in Fig. \ref{fig:mori2}.

\subsection{Debye-Hukel function}

The correlation function of the Gaussian polymer \cite{gross}, is given in
the closed form as
\begin{equation}
f(x)=\frac{2}{x}-\frac{2 (1-\exp (-x))}{x^2} \; .
\end{equation}

\begin{figure}
\begin{center}
\includegraphics[width=10cm]{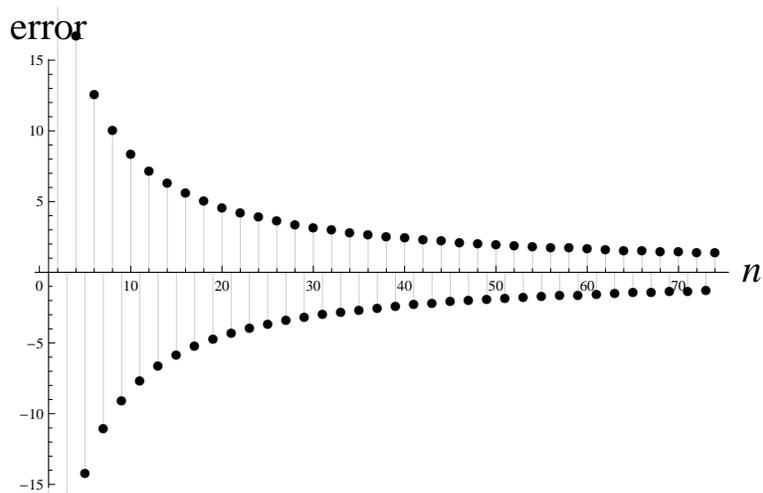}
\end{center}
\caption{
Debye-Hukel function. Percentage errors for the amplitude, with increasing 
approximation order, obtained with the standard Pad\'{e} approximants.
\label{fig:dh1}}.
\end{figure}

For small $x$, we have
\begin{equation}
f(x)\simeq 1 - \frac{x}{3} +\frac{x^2}{12} -\frac{x^3}{60} + \frac{x^4}{360} +O(x^5) \; ,
\end{equation}
and for large $x$, 
$$
f(x)\simeq 2 x^{-1} \qquad (x \ra \infty) \; .
$$
The standard Pad\'{e} scheme gives oscillating solutions with the magnitude
of oscillations slowly decreasing with the approximation order, as is shown in
Fig. \ref{fig:dh1}.

\begin{figure}
\begin{center}
\includegraphics[width=10cm]{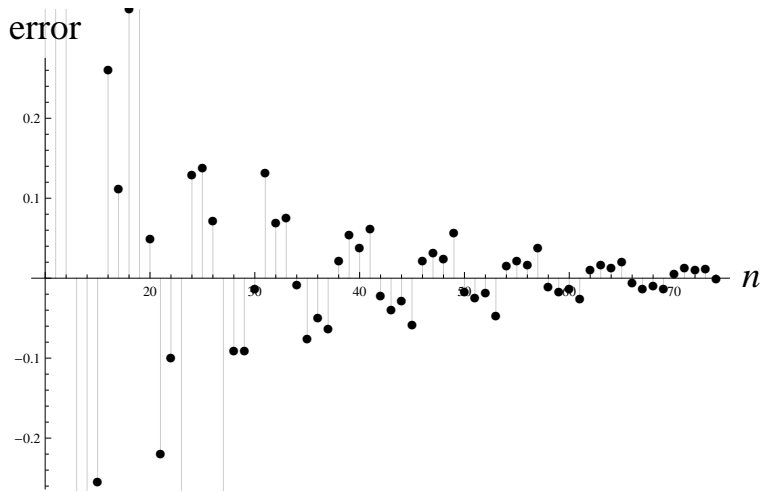}
\end{center}
\caption{
Debye-Hukel function. Percentage errors for the amplitude, with increasing 
approximation order, obtained with the corrected approximants.
\label{fig:dh2}}
\end{figure}

The method of corrected Pad\'{e} approximants, with a rather simple control function,
\begin{equation}
K(x)=R_{2}^{*}(x)=\frac{\sqrt{6}}{\sqrt{x^2+4 x+6}} \; ,
\end{equation}
demonstrates very good results improving convergence, as is shown 
in Fig. \ref{fig:dh2}.

\subsection{Branched polymer}

\begin{figure}
\begin{center}
\includegraphics[width=10cm]{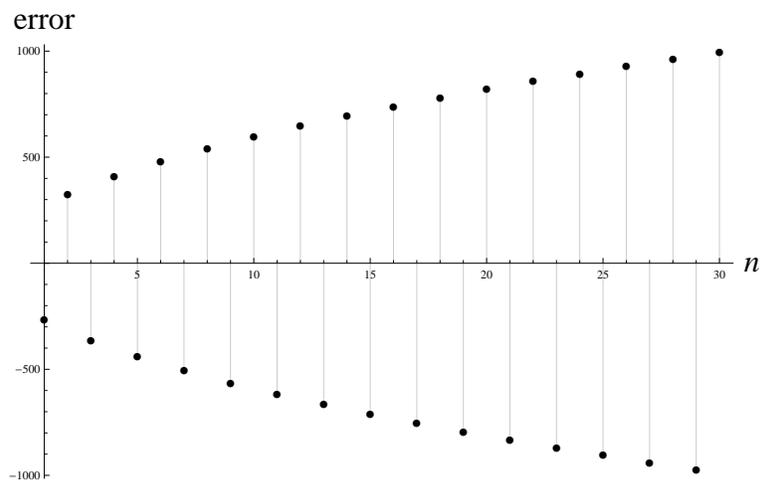}
\end{center}
\caption{
Branched polymer. Percentage errors for the amplitude, with increasing 
approximation order, obtained with Pad\'{e} approximants.
\label{fig:bran1}}
\end{figure}

The structure factor of a three-dimensional branched polymer is given by the
confluent hypergeometric function \cite{bran,bran1}
\begin{equation}
S(x)=F_1(1; \frac{3}{2}; \frac{3}{2} x) \; ,
\end{equation}
where $x$ is a dimensionless wave-vector modulus. At small $x$ (long-waves)
the coefficients in the expansion are given by the general expression
\begin{equation}
a_n= \left( -\frac{3}{2} \right)^n \Gamma (n+1)\; {\displaystyle /}\;
\frac{ \Gamma (1) \left[ n! \Gamma \left( n+\frac{3}{2}\right) \right ] }
{\Gamma \left(\frac{3}{2}\right) }  \; ,
\end{equation}
while for large $x$ (short-waves), one has 
$$
S(x)\simeq \frac{1}{3}\; x^{-1} \; .
$$
The standard Pad\'{e} scheme gives divergent oscillating results for the amplitude,
with an increasing magnitude for larger orders, as demonstrated in
Fig. \ref{fig:bran1}.

\begin{figure}
\begin{center}
\includegraphics[width=10cm]{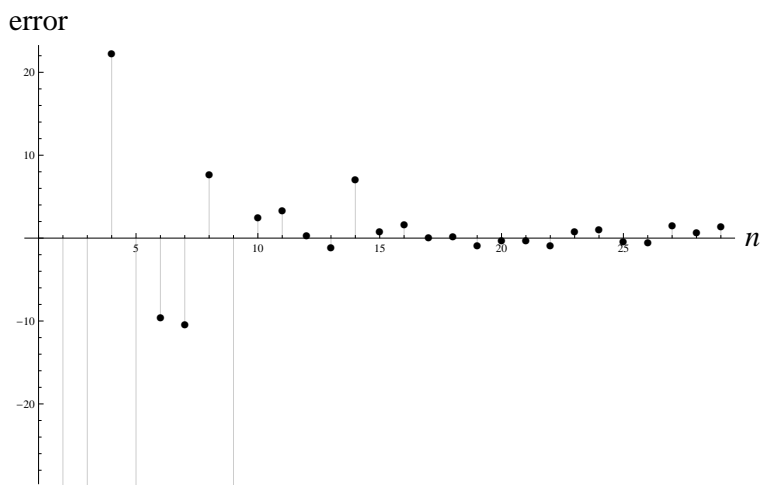}
\end{center}
\caption{
Branched polymer. Percentage errors for the amplitude, with increasing 
approximation order, obtained with the corrected Pad\'{e} approximants.
\label{fig:bran2}}
\end{figure}

The control function $K(x)$ in this case can be taken as the factor approximant
$$
K(x) = f_{3}^{*}(x)=
(1+(0.142857\, -0.255551 i) x)^{-0.5-1.67705 i} \times
$$
\begin{equation}
\times (1+(0.142857\, +0.255551 i) x)^{-0.5+1.67705 i} \; .
\end{equation}
The corrected Pad\'{e} approximants demonstrate very good convergence, as is
seen in  Fig. \ref{fig:bran2}.

\section{Treating convergence to wrong values}

It may happen that the standard Pad\'{e} method seems to converge very fast, 
however to a wrong value. This case can be viewed as divergent by oscillations, 
where the oscillations occur between this wrong value and $\infty$.

\subsection{Particle in a box}

In the calculation of the ground-state energy of a quantum particle in a one-dimensional 
box \cite{part}, one meets the following function 
\begin{equation}
\label{25}f(g)=\frac{\pi ^2}{128}\left(\frac{1}{2} +\frac{16}{\pi ^4g^2}+
\frac{1}{2} \sqrt{1+\frac{64}{\pi ^4g^2}}\right) \; .
\end{equation}
At small $g\rightarrow 0$, this function possesses the expansion
\begin{equation}
f(g)\simeq \frac 1{8\pi ^2g^2}\sum_{n=0}a_ng^n \; ,
\end{equation}
with 
$$
a_0=1, \quad a_1=\frac{\pi^2}4, \quad a_2=\frac{\pi^4}{32}, \quad a_{3}=\frac{\pi^6}{512}, 
\quad a_{4}=0, \quad a _{5}=-\frac{\pi^{10}}{131072}, \quad a_{6}=0, \ldots
$$
We shall be interested in finding the limiting value $f(\infty )$ in the
asymptotic expression
\begin{equation}
\label{26}f(g)\simeq f(\infty )\left( 1+O\left( g^{-2}\right) \right) ,\quad
\left( g\rightarrow \infty \right) \; ,
\end{equation}
whose exact limit is $0.077106$.

For convenience, instead of the original function (\ref{25}), we consider the
expression $8 \pi^2 g^2 \times f(g)$. The standard  Pad\'{e} scheme fails,
converging very fast, already at $n=2$, to the wrong value $0.0385531$, or 
producing rapid oscillations between this values $0.0385531$ and $\infty$, 
with varying the approximation order.

\begin{figure}
\begin{center}
\includegraphics[width=10cm]{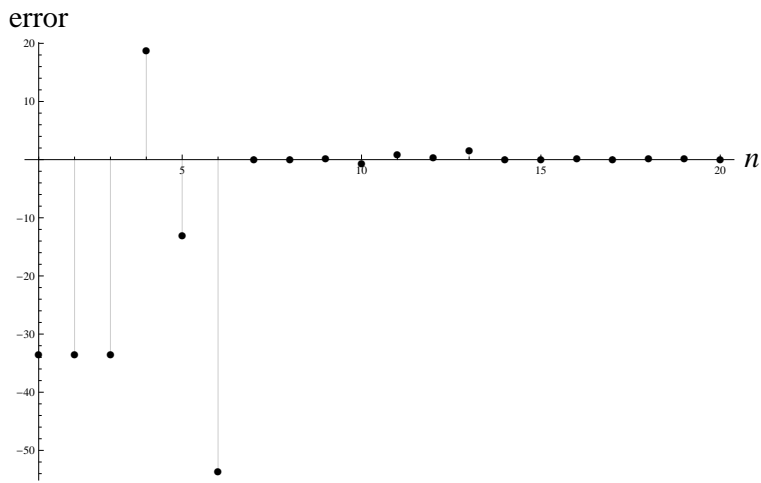}
\end{center}
\caption{
Particle in a box. Percentage errors for the amplitude, with increasing 
approximation order, for corrected Pad\'{e} approximants.
\label{fig:part}}
\end{figure}

The control function for the method of corrected Pad\'{e} approximants can 
be taken as the third-order iterated root approximant
\begin{equation}
K(g)=R_{3}^{*}(g)=\left(\frac{3 \pi ^6 g^3}{1024}+\left(\frac{\pi ^4 g^2}{64}+
\left(\frac{\pi ^2 g}{8}+1\right)^2\right)^{3/2}\right)^{2/3} \; .
\end{equation}
A very good convergence is achieved in the method of corrected Pad\'{e}
approximants, as is shown in Fig. \ref{fig:part}.

\subsection{Generating function}

The following function
\begin{equation}
\label{fgen}
f(x) = \left(\sqrt{x^2+1}+x\right)^a \; ,
\end{equation}
with an arbitrary $a$, is used as a touch-stone generating function \cite{genf}
for the sequence of numbers $a_n$ corresponding to the coefficients in the
expansion of (\ref{fgen}) at small $x$, with
\begin{equation}
\label{rec}
a_n=\frac{2^n \left(\frac{a}{2}-\frac{n}{2}+1\right)^{\bar{n}}}{n! \left(\frac{n}{a}+1\right)} \; ,
\end{equation}
while for large $x$, 
$$
f(x)\simeq 2^{a} x^{a} \qquad (x \ra \infty) \; . 
$$
Here $m^{\bar{k}}$ means $m (m+1)\ldots (k+m-1)$.

For concreteness, let us consider the case of $a=1/3$. Note that the function 
needed for constructing the standard scheme, that is, 
$$
T(x)\simeq 1-x+\frac{x^2}{2}-\frac{x^4}{8}+O(x^6) \; ,
$$
misses the powers of $x$ required for getting the correct large-variable limit. 
Therefore the value of the amplitude, predicted already at the first step $A_{1}=1$, 
does not change with $n$. Formally, the calculated amplitude oscillates between $1$ 
and $\infty$.

\begin{figure}
\begin{center}
\includegraphics[width=10cm]{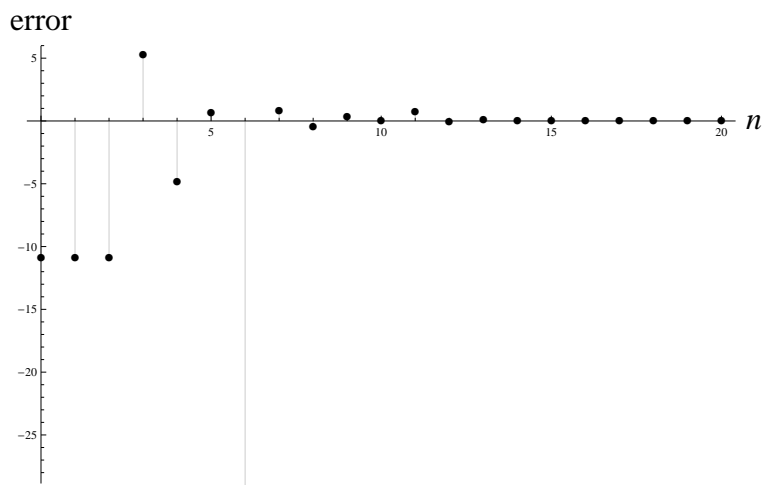}
\end{center}
\caption{
Generating function. Percentage errors for the amplitude, with increasing 
approximation order, obtained with the corrected approximants, corresponding 
to the control function (\ref{contgen}).
\label{fig:genf} }
\end{figure}

On the contrary, the corrected scheme works well, with the control function
\begin{equation}
\label{contgen}
K(x)=R_{2}^{*}(x)=\sqrt[6]{x^2+(x+1)^2} \; ,
\end{equation}
as is illustrated in Fig. \ref{fig:genf}.

\section{Treating strongly divergent solutions}

\subsection{Hard-core scattering problem}

Let us illustrate the method for the problem considered by Baker and
Gammel \cite{Baker_3}. When calculating the scattering length of a
repulsive square-well potential, one meets the integral
\begin{equation}
\label{sca}
 S(x) = \int_0^x \left ( \frac{\sin t}{t^3} \; - \; \frac{\cos t}{t^2}
\right ) \; dt \;  ,
\end{equation}
whose limit $S(\infty)=\pi/15$. Baker and Gammel state that the Taylor series 
expansion of this integral cannot be treated by the standard  Pad\'{e} techniques,
but require some more general methods. We show below that this integral can be 
effectively treated by means of the corrected approximants using the self-similar 
root approximants.

The small-variable expansion of this integral reads as
\begin{equation}
S(x) \simeq \frac{x}{9} -\frac{x^3}{135} + \frac{x^5}{2625}-
\frac{4x^7}{297675} + \frac{2x^9}{5893965} +O(x^{11}) \; .
\end{equation}

\begin{figure}
\begin{center}
\includegraphics[width=10cm]{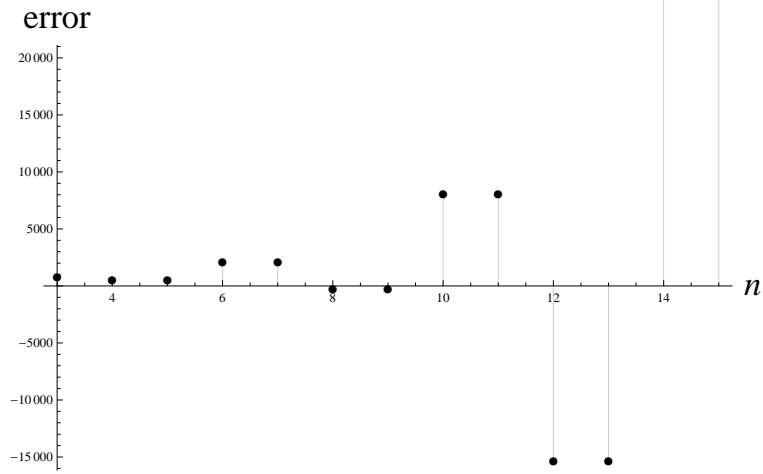}
\end{center}
\caption{
Scattering problem. Percentage errors for the amplitude, with increasing 
approximation order, obtained with the standard Pad\'{e} approximants.
\label{fig:baker}}
\end{figure}

The results of the standard Pad\'{e} scheme are meaninglessly divergent,
as is illustrated in Fig.\ref{fig:baker}.

It turns out that the sequence of the iterated roots gives the most stable results. 
The technique of iterated roots is very easy to apply, which yields \cite{g1,g2} 
the following few first terms,
\begin{eqnarray}
R_{1}^{*}(x)=\frac{x}{9 \sqrt{\frac{2 x^2}{15}+1}} \; ,
\nonumber
\\
R_{2}^{*}(x)=\frac{x}{9 \sqrt[4]{\frac{34 x^4}{2625}+\left(\frac{2 x^2}{15}+1\right)^2}} \; ,
\nonumber
\\
R_{3}^{*}(x)=\frac{x}{9 \sqrt[6]{\frac{152 x^6}{55125}+\left(\frac{34 x^4}{2625}+
\left(\frac{2 x^2}{15}+1\right)^2\right)^{3/2}}} \; ,...,
\end{eqnarray}
which can be extended to an arbitrary order.

Employing the corrected scheme, with the simplest root approximant $R^*_1(x)$ 
as a control function, we get the following sequence of approximations
for $S(\infty)$:
$$
S_1=0.30429, \quad S_2=0.247712, \quad S_3=0.238538, \quad S_4=0.238538, \quad
S_5=0.232624, 
$$
$$ 
S_6=0.228707, \quad S_7=0.225813, \quad S_8=0.223642, \quad 
S_9=0.221929, \quad S_{10}=0.220562, 
$$
$$
S_{11}=0.219428, \quad S_{12}=0.218486,
\quad S_{13}=0.217682, \quad S_{14}=0.216994, \quad S_{15}=0.216394, 
$$
$$
S_{16}=0.21587, \quad S_{17}=0.215405, \quad S_{18}=0.214992, \quad S_{19}=0.214621,
\quad S_{20}=0.214287, 
$$
$$
 S_{21}=0.213984, \quad S_{22}=0.213709, \quad 
S_{23}=0.213457, \quad S_{24}=0.213226, \quad S_{25}=0.213013 \; .
$$
The error of the last approximant equals $1.70644\%$.

\subsection{Wilson Loop}

The $N = 4$ super Yang - Mills circular Wilson loop \cite{wilson} is given by
the following expression,
\begin{equation}
\label{loop}
f(y)=\frac{2 \exp \left(-\sqrt{y}\right) I_1\left(\sqrt{y}\right)}{\sqrt{y}} \; ,
\end{equation}
where $I_{1}$ is a modified Bessel function of the first kind. Let us set
$\sqrt{y}=x$. For small $x$, we have
\begin{equation}
 f(x)\simeq 1-x+\frac{5 x^2}{8}-\frac{7 x^3}{24}+\frac{7 x^4}{64}+O(x^5) \; ,
\end{equation}
and at large $x$, one has
$$
f(x)\simeq \sqrt{\frac{2}{\pi }} x^{-3/2} \qquad (x\rightarrow \infty) \; .
$$

\begin{figure}
\begin{center}
\includegraphics[width=10cm]{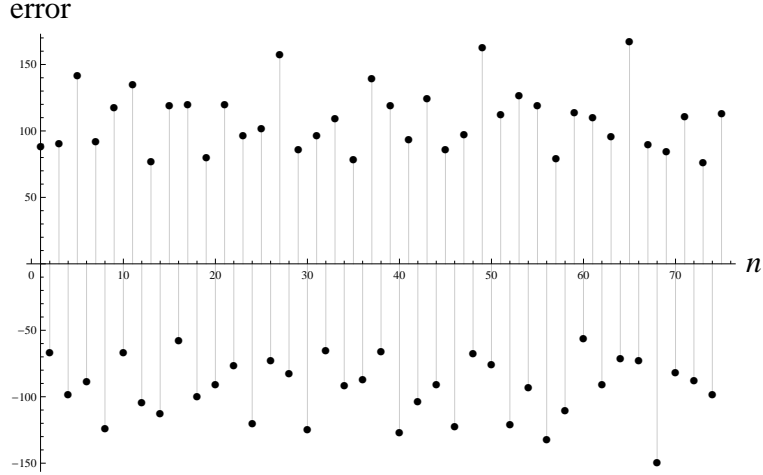}
\end{center}
\caption{
Wilson loop. Percentage errors for the amplitude, with increasing approximation 
order, obtained with the standard Pad\'{e} approximants.
\label{fig:loop1}}
\end{figure}

Applying the standard scheme leads to divergent results, as is clearly seen
in Fig. \ref{fig:loop1}.

\begin{figure}
\begin{center}
\includegraphics[width=10cm]{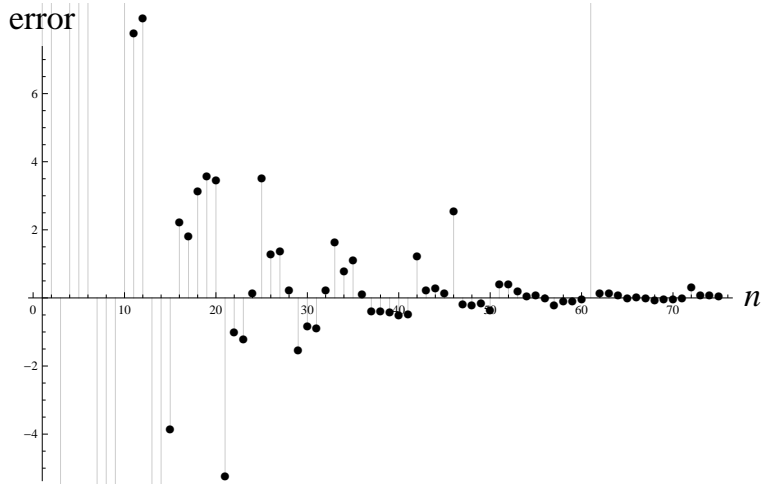}
\end{center}
\caption{
Wilson loop. Percentage errors for the amplitude, with increasing approximation 
order, obtained with the corrected approximants.
\label{fig:loop2}}
\end{figure}

Employing the control function
\begin{equation}
K(x)=R_{2}^{*}(x)=\left(\frac{5 x^2}{18}+\left(\frac{2 x}{3}+1\right)^2\right)^{-3/4} 
\end{equation}
and constructing the corrected Pad\'{e} approximants results in good numerical
convergence for the amplitude, as is shown in Fig. \ref{fig:loop2}.

\subsection{Error function}

The error function, often met in statistics and physics, is given by the integral,
\begin{equation}
F(x)=\int_0^x \exp \left(-u^2\right) \, du \; .
\end{equation}

\begin{figure}
\begin{center}
\includegraphics[width=10cm]{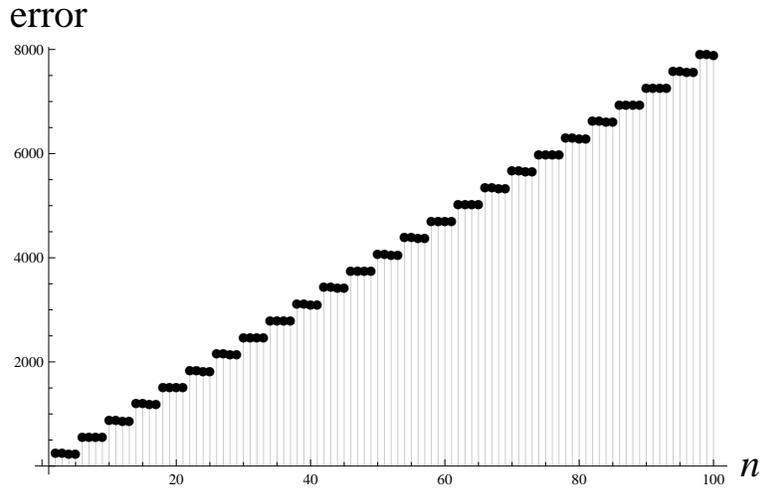}
\end{center}
\caption{
Error function. Percentage errors for the amplitude, with increasing approximation 
order, obtained with the standard Pad\'{e} approximants.
\label{fig:erf1}}
\end{figure}

For small $x$, we get
\begin{equation}
f(x)\simeq x-\frac{x^3}{3}+\frac{x^5}{10}+O(x^7) \; .
\end{equation}
And the limit at infinity is 
$$
F(\infty) =\frac{\sqrt{\pi }}{2} \; .
$$
In the framework of the standard Pad\'{e} scheme, one comes to a strongly
divergent sequence, as is seen in Fig. \ref{fig:erf1}.

\begin{figure}
\begin{center}
\includegraphics[width=10cm]{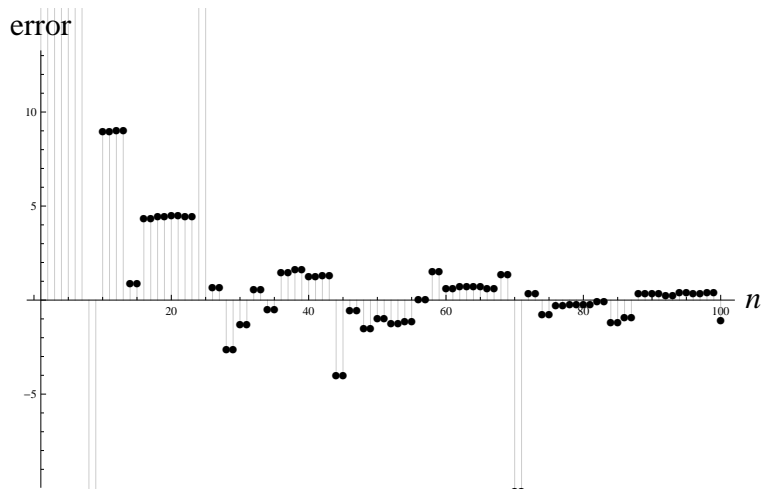}
\end{center}
\caption{
Error function. Percentage errors for the amplitude, with increasing approximation 
order, obtained with corrected approximants.
\label{fig:erf2}}
\end{figure}

In the framework of the corrected Pad\'{e} scheme, with the control function
\begin{equation}
K(x)=R_{3}^{*}(x)=\frac{x}{\sqrt[6]{\frac{16 x^6}{63}+\left(\frac{32 x^4}{45}+
\frac{4 x^2}{3}+1\right)^{3/2}}} \; ,
\end{equation}
good numerical convergence is achieved, as is shown in Fig. \ref{fig:erf2}.

\subsection{Debye function}

Let us consider the Debye function
\begin{equation}
D_{1}(x)=\frac{1}{x} \int_0^x \frac{y}{\exp (y)-1} \, dy \; .
\end{equation}
For large $x$, 
$$
D_{1}(x)\simeq A x^{-1} \; , 
$$
with the amplitude at infinity
$$
A=\frac{\pi ^2}{6}\approx 1.64493 \; .
$$
The expansion for small $x$ reads as
\begin{equation}
D_1(x)\approx 1-\frac{x}{4}+\frac{x^2}{36}-\frac{x^4}{3600}+\frac{x^6}{211680}+O(x^8) \; .
\end{equation}

\begin{figure}
\begin{center}
\includegraphics[width=10cm]{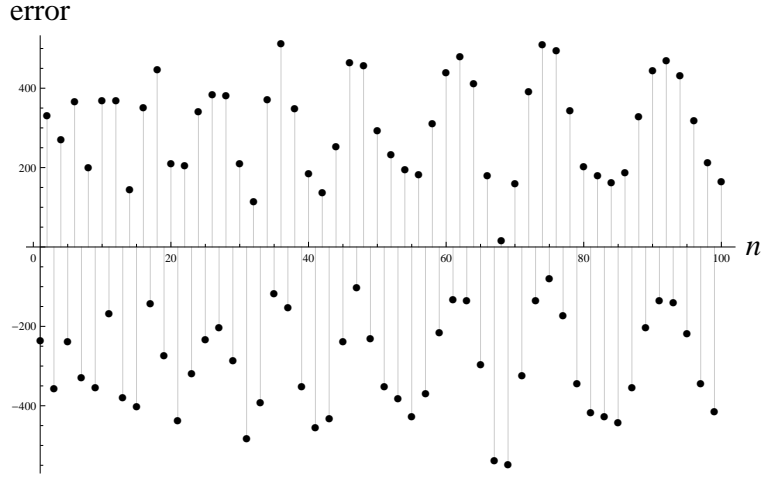}
\end{center}
\caption{
Debye function. Percentage errors for the amplitude, with increasing approximation 
order, obtained with the standard Pad\'{e} approximants.
\label{fig:debye1}}
\end{figure}

The standard Pad\'{e} approximations fail completely, following a rather
chaotic pattern, as is shown in Fig. \ref{fig:debye1}

\begin{figure}
\begin{center}
\includegraphics[width=10cm]{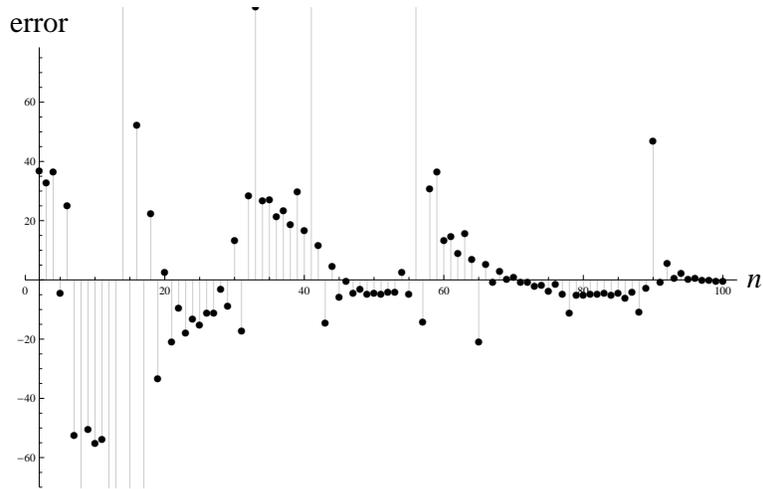}
\end{center}
\caption{
Debye function. Percentage errors for the amplitude, with increasing approximation 
order, obtained with the corrected Pad\'{e} approximants. The number $n$ changes 
from $n = 1$ to $n = 100$.
\label{fig:debye2}}
\end{figure}

The control function, playing the role of the initial approximation, is
\begin{equation}
K(x)=R_{2}^{*}(x)=\frac{1}{\sqrt{\frac{5 x^2}{72}+\left(\frac{x}{4}+1\right)^2}} \; .
\end{equation}
The corresponding corrected Pad\'{e} approximants work well, providing rather
good accuracy, as is shown in Fig. \ref{fig:debye2}.

\subsection{Connected moments}

An example from Ref. \cite{harm} concerns the application of the method of
connected moments to the calculation of the ground state energy of a harmonic
oscillator as the limit of generating function moments, as  "time" $t$ goes 
to infinity.

The generating function for the harmonic oscillator is
\begin{equation}
\label{con}
E(t)=\frac{121 u(t)^3+189199 u(t)^2+8180919 u(t)+6561}{(81-u(t)) \left(121 u(t)^2+20198 u(t)+
81\right)} \; ,
\end{equation}
where $u(t)=\exp (-4 t)$.

The $t$-expansion was also applied earlier to the case of $1d$ -antiferromagnet
\cite{antif}, but there are just a few starting terms in this $t$-expansion.

\begin{figure}
\begin{center}
\includegraphics[width=10cm]{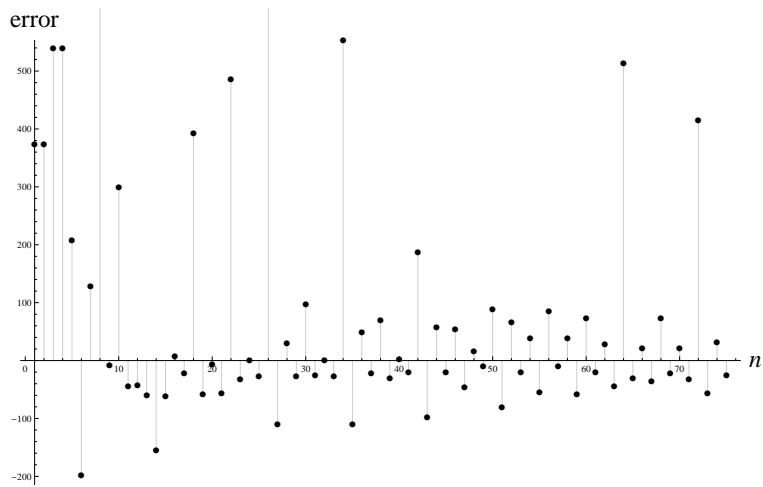}
\end{center}
\caption{
Generating function. Percentage errors for the amplitude, with increasing 
approximation order, obtained with the standard Pad\'{e} approximants.
\label{fig:harm1}}
\end{figure}

The standard Pad\'{e} scheme for eq. (\ref{con}) shows no convergence, as is evident
from the corresponding Fig. \ref{fig:harm1}.

\begin{figure}
\begin{center}
\includegraphics[width=10cm]{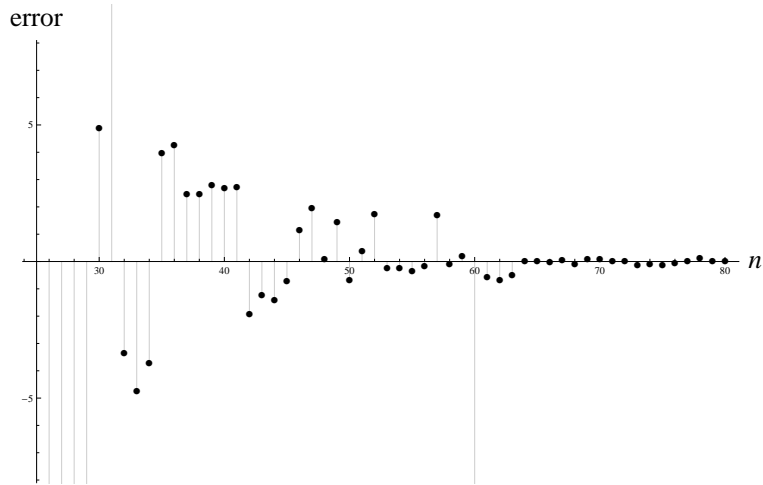}
\end{center}
\caption{
Generating function. Percentage errors for the amplitude, with increasing 
approximation order, obtained with the corrected Pad\'{e} approximants.
\label{fig:harm2}}
\end{figure}

Bur the corrected Pad\'{e} approximants demonstrate rather good convergence, with 
a simple control function corresponding to the shifted root approximant
\begin{equation}
K(t)=v_{2}+v_3 (v_4 t+1)^{-c} \; ,
\end{equation}
where 
$$
v_2=\frac{403171240048919}{85626857995920}, \quad v_3=\frac{36337990380139}{85626857995920},
\quad v_4=\frac{2331886111}{1340069829}, \quad c=\frac{9}{10} \; .
$$
The results of calculations, corresponding to the corrected Pad\'{e} scheme,
are shown in Fig. \ref{fig:harm2}.

\section{Massive Schwinger model in lattice theory}

In the previous sections, we have illustrated the method of corrected Pad\'{e} 
approximants for several simple cases. Now we show its applicability for more
complicated cases that are of physical interest.  

In the present section, we consider the massive Schwinger model in Hamiltonian
lattice theory \cite{Schwinger_41,Hamer_42}. This model describes quantum 
electrodynamics in two space-time dimensions. Its fascinating features include
many of the properties of quantum chromodynamics, such as confinement, chiral
symmetry breaking, and a topological vacuum. Because of these effects, the 
Schwinger model has attracted great interest. 

Let us consider the energy gap between the lowest and first excited states of the
vector boson, which can be represented as a function $\Delta(z)$ of the variable
\be
\label{M1} 
 z = x^2 \; , \qquad x \equiv \frac{1}{g^2 a^2} \;  ,
\ee
where $g$ is a coupling parameter and $a$ is lattice spacing. The gap can be 
expressed as an expansion 
\be
\label{M2}
\Dlt(z) \simeq \sum_n a_n z^n \qquad (z \ra 0 )
\ee
in powers of asymptotically small $z$, when the coupling parameter is strong. 
The coefficients here are
$$
a_0 = 1 \; , \qquad a_1 = 2 \; , \qquad a_2 = -10\; , 
$$
$$
a_3 = 78.66667 \; , \qquad a_4 = -736.2222 \; , \qquad a_5 = 7572.929\; , 
$$
$$
a_6 = -82736.69 \; , \qquad a_7 = 942803.4 \;   .
$$
The strong increase of the coefficients makes the series in powers of $z$ 
widely divergent. 

However, the transition from the lattice formulation to the continuous limit
requires taking the limit $a \ra 0$, which implies $z \ra \infty$. In this limit,
the gap behaves as
\be
\label{M3} 
 \Dlt(z) \simeq 0.5642 z^{1/4} \qquad ( z \ra \infty) \;  .
\ee

Using Pad\'{e} approximants, we get $A_7 = 0.680043$, with the percentage error
of $21 \%$. While employing the corrected Pad\'{e} approximants, with the control 
function
$$
K(z) = R_2^*(z) = \left (  ( 1 + 8z )^2 - 32 z^2 \right )^{1/8} \;   ,
$$
we find $A_7 = 0.591181$, with the error of $4.8 \%$, which is essentially better
than the standard approach.

\section{Critical temperature of weakly interacting Bose gas}

Three-dimensional Bose gas, at low temperature experiences Bose-Einstein 
condensation at a critical temperature $T_c$. The properties of Bose-condensed 
gas have recently been intensively studied both theoretically and experimentally 
(see recent books and review articles 
(\cite{Pethick_41,Lieb_42,Letokhov_43,Y_44,Y_45,Y_46}). An interesting problem 
is the calculation of the critical temperature $T_c$ of weakly interacting uniform 
Bose gas, as compared to the Bose-condensation temperature of ideal gas
\be
\label{C1}
 T_0 = \frac{2\pi}{m} \; \left [ \frac{\rho}{\zeta(3/2)} \right ]^{2/3} \; ,
\ee
where $m$ is atomic mass and $\rho$, average density. The Planck and Boltzmann 
constants are set to one.

The difficulty of this problem is that the ideal Bose gas and interacting Bose gas 
are in different classes of universality, so that the transition from the ideal 
gas to the interacting one cannot be done by simple perturbation theory \cite{Holzmann_47}.

One considers the relative critical temperature shift
\be
\label{C2}
 \frac{\Dlt T_c}{T_0}  \equiv \frac{T_c - T_0}{T_0} \;   ,
\ee
which, at weak interaction, can be presented by the expression
\be
\label{C3}
 \frac{\Dlt T_c}{T_0} \simeq c_1 \gm \qquad ( \gm \ra 0 ) \;  ,
\ee
as a function of the gas parameter
\be
\label{C4}
 \gm \equiv \rho^{1/3} a_s \;  ,
\ee
with $a_s$ being scattering length. 

The coefficient $c_1$, for the two-component $O(2)$ field theory, describing 
the real interacting Bose gas, has been calculated by different ways, including 
Ursel operator techniques \cite{Baym_48} ($c_1 = 1.95$), renormalization-group 
methods \cite{Ledowski_49,Blaizot_50,Benitez_51} ($c_1 = 1.15 - 1.37$). Optimized 
perturbation theory \cite{Y_52,Y_53} has been used for calculating $c_1$ by 
either introducing control functions into an initial Lagrangian 
\cite{Cruz_54,Cruz_55,Cruz_56,Kneur_57,Kneur_58,Kneur_59} or into a variable
transformation \cite{Kleinert_60,Kastening_61,Kastening_62,Kastening_63}. 
Control functions are defined by optimization conditions, among which the most 
often used are the variational and finite-difference conditions \cite{Y_64,Y_65,Y_66}. 
These conditions are equivalent to each other and, for defining $c_1$, lead to 
very close results \cite{Kneur_57}. Optimized perturbation theory, with increasing 
order, numerically converges \cite{Kneur_57,Braaten_67,Braaten_68} to the value 
of the coefficient $c_1$ close to that found in Monte Carlo simulations 
\cite{Arnold_69,Arnold_70,Arnold_71} ($c_1 = 1.32 \pm 0.02$) and 
\cite{Kashurnikov_72} ($c_1 = 1.29 \pm 0.05$). The critical temperature for larger 
values of the gas parameter has also been studied by Monte Carlo simulations 
\cite{Pilati_73}. And the critical temperature of the $O(N_1) \times O(N_2)$ 
field theory has been calculated by optimized perturbation theory \cite{Pinto_74}. 
A more detailed account of different attempts of calculating $T_c$ can be found in 
Refs. \cite{Kleinert_75,Andersen_76,Y_77}. 

This brief discussion illustrates the importance of calculating the critical 
temperature in the $O(N)$ field theory. Here we show how the critical temperature can 
be easily found by the method of corrected Pad\'{e} approximants.

Let us start with the seven-loop expansion \cite{Kastening_63} for the 
parameter $c = c(x)$,
\be
\label{C5}
 c_1(x) \simeq a_1 x + a_2 x^2 + a_3 x^3 + a_4 x^4 + a_5 x^5 \;  ,
\ee
derived for an asymptotically small variable $x \ra 0$, with
$$
a_1 = 0.223286 \; , \qquad a_2 = -0.0661032 \; , \qquad
a_3 = 0.026446 \; ,
$$
$$
a_4 = -0.0129177 \; , \qquad a_5 = 0.007290373 \;   .
$$
But the sought value of $c_1$ corresponds to the limit
\be
\label{C6}
 c_1 = \lim_{x\ra\infty} c_1(x) \;  .
\ee
  
In the calculations of Kastening \cite{Kastening_63}, the variable transformation 
$x = z/ (1-z)^{1/\omega}$ is used, so that $x \ra \infty$ as $z \ra 1$. The 
obtained results essentially depend on the chosen value of $\omega$, for different 
choices varying between $c_1 = 1.161$ and $c_1 = 1.376$. It is also worth mentioning
that it would be possible to invoke other variable transformations. However, 
as is known \cite{Bailey_78}, different variable changes can strongly influence
the sought limits.      

The best Pad\'{e} approximant $P_{2/3}$ yields the value $c_1 = 0.982$,
which is rather far from the Monte Carlo results. Employing the method 
of corrected Pad\'{e} approximants, with the control function 
$$
K(x) = R_2^*(x) = \frac{0.223286x}{((1+0.296x)^2-0.0616x^2)^{1/2}} \;   ,
$$
we find $c_1 = 1.386$, which is closer to the Monte Carlo results.  

In the same way, it is straightforward to find the values of $c_1$ for the
$O(1)$ field theory, for which the seven-loop expansion \cite{Kastening_63} 
gives form (\ref{C5}) with the coefficients
$$
a_1 = 0.334931 \; , \qquad a_2 = -0.178478 \; , \qquad
a_3 = 0.129786 \;   ,
$$
$$
a_4 = -0.115999 \; , \qquad a_5 = 0.120433 \; .
$$
The best Pad\'{e} approximant gives $c_1 = 0.824$, which is much lower than
the Monte Carlo result \cite{Sun_79} equal to $c_1 = 1.09 \pm 0.09$. Resorting
to the corrected Pad\'{e} approximants, with the control function $R_2^*(x)$,
we find $c_1 = 1.207$, which is a bit higher than the Monte Carlo value. 

For the $O(4)$ field theory, the seven-loop expansion yields the coefficients
$$
a_1 = 0.167465 \; , \qquad a_2 = -0.0297465 \; , \qquad
a_3 = 0.00700448 \; ,
$$
$$
 a_4 = -0.00198926 \; , \qquad a_5 = 0.000647007 \;  .
$$
The best Pad\'{e} approximant yields $c_1 = 1.219$, which is lower than the 
Monte Carlo result \cite{Sun_79} giving $c_1 = 1.6 \pm 0.1$. By using the 
method of corrected Pad\'{e} approximants, again with the control function
$R_2^*(x)$, we obtain $c_1 = 1.6$ coinciding with the Monte Carlo result.

\section{Equation of state for hard-sphere fluids}

The structure of many real fluids is mainly determined by repulsive forces,
because of which the hard-sphere fluid is the simplest and the most widely 
used model for describing the behaviour of real fluids. 

The equation of state of hard-sphere fluids is characterized \cite{Mulero_80}
by the so-called compressibility factor
\be
\label{E1}
Z \equiv \frac{P}{\rho k_B T} = Z(y) \qquad 
\left ( y \equiv \frac{\pi\rho}{6} \right ) \;   ,
\ee
where $P$ is pressure, $\rho$ is density, $T$ is temperature, and 
$y \equiv \pi \rho/6$ is packing fraction. This factor is usually represented
by the virial expansion
\be
\label{E2}
 Z = 1 + \sum_{n=2}^\infty B_n y^{n-1} \;  .
\ee
The virial coefficients $B_n$ do not depend on temperature and are defined in 
terms of integrals whose integrands are products of Mayer functions. Only the 
first four virial coefficients can be calculated analytically:
$$
B_1 = 1 \; , \qquad B_2 = 4 \; , \qquad B_3 = 10 \; , 
$$
$$
B_4 = \frac{2707\pi+[438\sqrt{2} - 4131\arccos(1/3)]}{70\pi} \; =
18.364768 \;   .
$$
The higher virial coefficients have been calculated \cite{Clisby_81,Clisby_82} 
numerically:
$$
B_5=28.224512 \; , \qquad B_6 = 39.815148 \; , \qquad B_7 = 53.344420 \; ,
$$
$$
B_8=68.537549 \; , \qquad B_9 = 85.812838 \; , \qquad B_{10} = 105.775104 \; ,
$$
$$
B_{11} = 127.93 \; , \qquad  B_{12} = 152.67 \; , \qquad B_{13} = 181.19 \; ,
$$
$$
B_{14} = 214.75 \; , \qquad  B_{15} = 246.96 \; , \qquad B_{16} = 279.17 \;   .
$$
The values of the virial coefficients increase with the expansion order,
because of which this expansion is strongly divergent. 

In order to get an equation of state, one either employs Pad\'{e} 
approximants complimented by some phenomenological terms or constructs purely 
phenomenological equations. The standard way of checking the validity of the 
so constructed equations is as follows. One forms an equation exactly reproducing 
the first ten virial coefficients. Then this equation is expanded in powers 
of the packing fraction $y$ and one examines how such an expansion 
reproduces the last virial coefficients from $B_{11}$ to $B_{16}$. The maximal error 
in reproducing these last virial coefficients defines the accuracy of the 
studied equation of state. The analysis of a great number of different equations 
of states has been accomplished in Refs. \cite{Wu_83,Tian_84,Tian_85}. Referring 
to these papers, we give below a brief account of their findings. 

One of the most popular equations is the Cannahan-Starling equation, which has 
the structure of a $P_{3/3}$ Pad\'{e} approximant. This equation reproduces 
the higher virial coefficients with a maximal error of $3.6 \%$. 

Clisby and McCoy proposed the $P_{4/5}$ and $P_{5/4}$ Pad\'{e} approximants,
whose maximal errors in the reproduction of higher virial coefficients are 
$9.1 \%$ and $7.1 \%$, respectively.

Kolafa, Labik, and Malijevsky suggested the equations in the form of series in
powers of $y/(1-y)$, with the parameters fitted to the lower densities. Two 
fitting procedures lead to the reproduction of the higher virial coefficients
with the errors $75 \%$ and $34 \%$. 

Liu suggested a phenomenological combination of Pad\'{e} approximants and
power series, with the maximal error of $15 \%$. 

Santos and Lopez de Haro proposed a phenomenological equation, having the 
error of $8.9 \%$.     

Tian, Jiang, Gui, and Mulero \cite{Tian_84} constructed a phenomenological 
equation as a Laurent series in powers of $y-b$. As they write, it was possible
to construct $57$ different variants of such equations, the best of which 
reproduces the higher virial coefficients with an error of $3.2 \%$. 

Applying the method of corrected Pad\'{e} approximants, it is convenient, first,
to introduce the notation
\be
\label{E3}
 x = \frac{y}{1-y} \qquad \left ( y = \frac{x}{1+x} \right ) \;  ,
\ee
such that $x \ra \infty$, when $y \ra 1$. This allows us to use the standard 
techniques of Sec. 2 by considering the limit $x \ra \infty$. 

As follows from some phenomenological equations \cite{Wu_83,Tian_84,Tian_85}, 
for instance, from the Cannahan-Starling equation, the compressibility factor   
in the vicinity of $y = 1$ behaves as
$$
 Z(y) \simeq \frac{2}{(1-y)^3} \qquad ( y \ra 1 ) \;  .
$$
Hence, in terms of the variable $x$, we have
\be 
\label{E4}
 Z(y(x)) \simeq 2 x^3 \qquad ( x\ra \infty) \;  .
\ee

Following the scheme of Sec. 2, we take as a control function the second-order
root approximant
\be
\label{E5}
 K(x) = R_2^*(x) = \left ( \left ( 1 + \frac{4}{3} \; x \right )^2 +
\frac{4}{9} \; x^2 \right )^{3/2} \;  .
\ee
For the ratio $Z(y)/K(x(y))$, we construct a corrected Pad\'{e} approximant
$P_{5/5}(y)$ taking into account the first ten virial coefficients. Then the 
resulting equation of state is
\be
\label{E6}
 Z^*(y) = K(x(y)) P_{5/5}(y) \;  .
\ee
This equation predicts the higher virial coefficients, up to $B_{16}$, with an 
error not exceeding $2.5 \%$. 

Note that the method of estimating the accuracy of the approach by predicting
the higher-order coefficients goes back to Feynman, as discussed in Ref. 
\cite{Samuel}.

\section{Fluctuating pressure of fluid membrane}

Membranes are frequent structures in various biological, chemical, and mechanical
systems, whose description is not easy. One of the often met types is a fluid 
membrane between two fixed boundaries. To calculate the pressure of this membrane,
one introduces a potential of strength $\lambda$ between two walls at distance $L$.
Then the membrane pressure, in dimensionless units, can be represented as a function 
$p(x)$ of a variable $x = 1/(\lambda L)^2$, yielding \cite{part} the expansion
\be
\label{F1}
 p(x) \simeq \frac{\pi^2}{8x^2} \sum_{n=0}^k a_n x^n \;  ,
\ee
when $x \ra 0$, with the coefficients
$$
a_0 = 1 \; , \qquad a_1 = \frac{1}{4} \; , \qquad a_2 = \frac{1}{32} \; , \qquad
a_3 = 2.176347\times 10^{-3} \; , 
$$ 
$$
 a_4 = 0.552721 \times 10^{-4} \; , \qquad  a_5 = -0.721482 \times 10^{-5} \; , 
\qquad  a_6 = -1.777848 \times 10^{-6} \; .
$$
But in order to return to the case of a membrane between two rigid walls, one has 
to remove the auxiliary potential, setting $\lambda \ra 0$, when $x \ra \infty$. 
Hence, one needs to find the limit
$$
 p = \lim_{x\ra\infty} p(x) \;  .
$$
Several attempts employing directly Pad{\'{e} approximants, as well as other 
approximants \cite{Y_88}, have been shown to give bad accuracy, with the best
Pad{\'{e} approximants giving $p = 0.015$, which is quite far form the value
$$
p = 0.0798 \pm 0.0003
$$
obtained in Monte Carlo simulations \cite{Gompper_90}. A closer to the Monte Carlo
result was found by Kastening \cite{part} by employing optimized perturbation theory, 
with introducing control functions through a change of variables \cite{Kleinert_75}. 
As has been mentioned above, the results essentially depend on the kind of the 
variable transformation. Kastening found $p = 0.0821$, which is $3 \%$ higher than 
the Monte Carlo value. Below, we show how an accurate value for the limit $p(\infty)$ 
can be easily derived by means of corrected Pad{\'{e} approximants.    

Following the method of corrected Pad{\'{e} approximants, with the control 
function in the form of a root approximant
\be
\label{F2}
 K(x) = \left ( \left ( \left ( 1 + \frac{1}{8} \; x \right )^2 + 
\frac{1}{64} \; x^2 \right )^{3/2} + A_3 x^3 \right )^{2/3} \;  ,
\ee
where $A_3 = 0.00326452$, we consider the series for the ratio $p(x)/K(x)$ and  
construct the Pad{\'{e} approximant $P_{4/4}$, so that the pressure becomes 
\be
\label{F3}
 p(x) = \frac{\pi^2}{8x^2} \; K(x) P_{4/4}(x) \;  .
\ee
Taking here the limit $x \ra \infty$, we obtain 
$$
p(\infty) = 0.0806 \;   ,
$$
which is very close to the Monte Carlo value, deviating from it only by $1 \%$.

\section{Conclusion}

We have suggested a method of corrected Pad\'{e} approximants, whose idea is 
based on using as an initial approximation a self-similar approximant taking
into account irrational features of the sought function. We showed that the
suggested method covers a much broader class of functions than the standard 
Pad\'{e} approximants do, including the cases of indeterminate problems. 
The corrected Pad\'{e} approximants are easy to construct. Their computational 
accuracy, when the standard Pad\'{e} scheme appears to converge, is either 
close to that of diagonal Pad\'{e} approximants, when these exist, or is 
essentially better, making the rate of convergence faster. 

But most importantly, the corrected scheme works well, even when the standard
Pad\'{e} scheme completely fails, being either divergent or not defined at all.

We have illustrated the advantage of the suggested corrected Pad\'{e}
approximants by several simple functions related to physical applications. Also,
we demonstrated the efficiency of the method for several complicated physical
problems, such as the massive Schwinger model in lattice field theory, the scalar
field theory $O(N)$ having to do with the phase transition in Bose gases, the 
equation of state for hard-sphere fluids, and fluctuating pressure of a fluid 
membrane. The method of corrected Pad\'{e} approximants is shown to combine rather
easy calculations with good accuracy.

\end{document}